\documentclass{article}%
\usepackage{amsmath}%
\usepackage{amsfonts}%
\usepackage{amssymb}%
\usepackage{graphicx}

\begin{document}

\title{On Methods for Transforming and Solving Finite Series}
\author{Henrik Stenlund\thanks{The author is obliged to Visilab Signal Technologies for supporting this work.}
\\Visilab Signal Technologies Oy, Finland}
\date{28th January, 2016}
\maketitle
\begin{abstract}
In this work we present new methods for transforming and solving finite series by using the Laplace transform. In addition we introduce both an alternative method based on the Fourier transform and a simplified approach. The latter allows a quick solution in some cases. \footnote{Visilab Report \#2016-01}
\subsection{Keywords}
Summation of finite series, infinite series, Laplace transform, inverse Laplace transform, Fourier transform
\subsection{Mathematical Classification}
MSC: 44A10, 11M41, 16W60, 20F14, 40A25, 65B10
\end{abstract}
\tableofcontents

\section{Introduction}
\subsection{General}
The motivation for this paper has been the need to handle a general finite series of the form below,
\begin{equation}
A_N=\sum_{k=1}^{N}{g(k)} \label{eqn2}
\end{equation}
It needs to be solved or transformed into another form which is easier in further processing. Just a few methods exist of general use for handling finite series. The most important is the Euler-Maclaurin formula. Even a transform of the finite series may become useful as an intermediate step in solving it. The new form can be another finite or infinite series or an integral. There exists only a relatively small number of tabularized finite series. 
\subsection{The Euler-Maclaurin Formula}
The Euler-Maclaurin formula \cite{Abramowitz1970}, \cite{Jeffrey2008}, \cite{Ivic1985} is the most important method for solving finite series. It presents an expression between the finite series and integrals as follows, assuming $a,b$ to be integers. Let the first $2n$ derivatives of $f(x)$ be continuous on an interval $[a,b]$. The interval is divided into equal parts $h=\frac{(a-b)}{n}$. Then for some $\theta , 0\leq{\theta}\leq{1}$ 
\begin{equation}
\sum_{k=0}^{m}{f(a+kh)}=\frac{1}{h}\int^{b}_{a}{f(x)dx}+\frac{1}{2}(f(a)+f(b))+\sum_{k=1}^{n-1}{B_{2k}h^{2k-1}\frac{[f^{(2k-1)}(b)-f^{(2k-1)}(a)]}{(2k)!}}+ \nonumber
\end{equation}
\begin{equation}
+\frac{h^{2n}}{(2n)!}{B_{2n}}\sum_{k=0}^{m-1}{f^{(2n)}(a+kh+\theta{h})} \label{eqn33}
\end{equation}
The $B_{i}$ are Bernoulli numbers. 
\subsection{A Trivial Summation Method}
A simple method for summation of finite series is the following. 
We have a difference 
\begin{equation}
\Delta{u(k)}=u(k+1)-u(k)=f(k) \label{eqn80}
\end{equation}
We can subject it to the summation operator
\begin{equation}
\Delta^{-1}{\Delta{u_k}}=\Delta^{-1}{f(k)} \label{eqn90}
\end{equation}
which is equal to 
\begin{equation}
\sum_{k=1}^{N}{f(k)}=u(N+1)-u(1)  \label{eqn94} 
\end{equation}
This is the sum of the difference of the $u(k)$. Traditional summation methods may be helpful in simple cases, found in \cite{Schaum1971}.

In Section 2 we derive the method, its generalization and explain the process of using the method. We exhibit also some illustrating sample cases. In Section 3 we show the Fourier method and Section 4 contains the simplified method. Appendix A shows a few extensions for the extended method better suit various cases. 

\section{The Laplace Transform Method for Finite Series}
\subsection{The Simple Form}
We follow the steps of a method for infinite series, recently developed by the author \cite{Stenlund20142}. The \textbf{infinite series method} has the following simple form
\begin{equation}
\sum_{k=1}^{\infty}{g(k)}=\int^{\infty}_{0}{\frac{dt{G(t)}}{e^{t}-1}} \label{eqn100}
\end{equation}
and the extended form is
\newpage
\begin{equation}
\int^{\infty}_{0}{\frac{dt{G(t)}}{e^{\alpha{t}}-1}}\  \  =\  \  \sum_{k=1}^{\infty}{g(\alpha{k})} \  \  =\  \  f(\alpha) \label{eqn130}
\end{equation}
\begin{equation}
       {\textsl{L}_{x},{\alpha}} \Uparrow \  \  \   \Downarrow{\textsl{L}^{-1}_{\alpha},{x}}   \nonumber
\end{equation}
\begin{equation}
\sum_{k=1}^{\infty}{\frac{G(\frac{x}{k})}{k}}\  \  =\  \  F(x) \label{eqn140}
\end{equation}
The method for \textbf{finite series} is developed analogously but has significant differences. It is important to note that we don't allow $N\rightarrow{\infty}$. $N$ can be even one and the index function $g(k)$ can be unusual. $g(k)$ must be a mapping of the index, not a set of random numbers. The function must be a surjection if not an injection or bijection. $g(k)$ can be of such a form that while letting $N\rightarrow{\infty}$ the resulting infinite series may be diverging. If the limit $N\rightarrow{\infty}$ is important then the infinite series method should be consulted instead. The most important requirement is that $g(k)$ must have an inverse Laplace transform $G(t)$. 
\begin{equation}
g(k)=\int^{\infty}_{0}{e^{-k{t}}{G(t)}dt}  \label{eqn354}
\end{equation}
Thus
\begin{equation}
\sum_{k=1}^{N}{g(k)}=\sum_{k=1}^{N}\int^{\infty}_{0}{e^{-k{t}}{G(t)}dt} \label{eqn356}
\end{equation}
Since we have a finite number of terms in the series in all cases, we can interchange the order of summation and integration
\begin{equation}
\sum_{k=1}^{N}{g(k)}=\int^{\infty}_{0}{\sum_{k=1}^{N}e^{-k{t}}{G(t)}dt}=\int^{\infty}_{0}{dt{G(t)}\sum_{k=1}^{N}e^{-k{t}}} \label{eqn358}
\end{equation}
We know that
\begin{equation}
\sum_{k=1}^{N}e^{-k{t}}=\frac{1-e^{-t{N}}}{e^{t}-1} \label{eqn362}
\end{equation}
which is very simple to prove and the finite series becomes
\begin{equation}
\sum_{k=1}^{N}{g(k)}=\int^{\infty}_{0}{\frac{dt{G(t)({1-e^{-t{N}}})}}{e^{t}-1}} \label{eqn364}
\end{equation}
We insist on the following requirements:
\begin{itemize}
	\item 1. $g(k)$ has an inverse Laplace transform $G(t)$
	\item 2. the resulting definite integral is finite
\end{itemize}
\subsection{Example Cases}
\subsubsection{A Simple Demonstration of the Method} 
The method is tested in the following by applying it to a finite series of a trivial delta function. 
\begin{equation}
\sum_{k=1}^{N}{\frac{\delta(a-\frac{x}{k})}{k}}  \label{eqn800}
\end{equation}
We solve the  finite series by using the equations ($\ref{eqn428}$) and ($\ref{eqn430}$). Thus 
\begin{equation}
G(\frac{x}{k})=\delta(a-\frac{x}{k})   \label{eqn802}
\end{equation}
We can immediately calculate the integral as
\begin{equation}
\int^{\infty}_{0}{\frac{dt{\delta(a-t)({1-e^{-\alpha{tN}}})}}{(e^{\alpha{t}}-1)}}=\frac{({1-e^{-\alpha{Na}}})}{(e^{a\alpha}-1)} \label{eqn812}
\end{equation}
\begin{equation}
=\sum_{k=1}^{N}{g(\alpha{k})}=f(\alpha) \nonumber
\end{equation}
On the other hand
\begin{equation}
g(k)=\int^{\infty}_{0}{dt{e^{-kt}G(t)}}=e^{-ka} \label{eqn850}
\end{equation}
Thus we are led to the old identity
\begin{equation}
f(\alpha)=\sum_{k=1}^{N}{e^{-ka\alpha}}=\frac{({1-e^{-\alpha{Na}}})}{(e^{a\alpha}-1)} \label{eqn860}
\end{equation}
By the inverse Laplace transform we can get $F(x)$ from this
\begin{equation}
F(x)=L^{-1}_{\alpha}[\sum_{n=0}^{\infty}{e^{-ka\alpha(n+1)}-e^{-ka\alpha(N+n+1)}}]=  \label{eqn870}
\end{equation}
\begin{equation}
=\sum_{n=0}^{\infty}{[\delta(x-a(n+1))-\delta(x-a(N+n+1))]}  \label{eqn871}
\end{equation}
The terms cancel each other when $n>N$ and we are left with, by using the properties of the delta function
\begin{equation}
F(x)=\sum_{n=1}^{N}{\frac{\delta(a-\frac{x}{n})}{n}}  \label{eqn880}
\end{equation}
The loop is closed. 
\subsubsection{A Fractional Finite Series} 
In the following we attempt to solve a fractional finite series by using the equations ($\ref{eqn410}$) and ($\ref{eqn415}$)
\begin{equation}
\sum_{k=1}^{N}{\frac{a}{(k^2+a^2)}}=\sum_{k=1}^{N}{g(k)} ,a\neq{0} \  \label{eqn1850}
\end{equation}
The inverse Laplace transform of $g(k)$ is
\begin{equation}
G(t)=sin(at) \label{eqn1852}
\end{equation}
We use our formulas to get with the parameter $\alpha$
\begin{equation}
\sum_{k=1}^{N}{g(\alpha{k})}=\int^{\infty}_{0}{\frac{dt{({1-e^{-\alpha{tN}}})sin(at)}}{(e^{\alpha{t}}-1)}} \label{eqn1855}
\end{equation}
Expanding the $sin$ function as exponential functions and then as power series in the integral, will give us after changing the variable to $x=\alpha{t}$
\begin{equation}
\sum_{k=1}^{N}{\frac{a}{(k^2+a^2)}}= \nonumber
\end{equation}
\begin{equation}
=\frac{1}{2i}\sum_{n=1}^{\infty}{\frac{{\alpha^{n-1}}{[(ia)^n-(-ia)^n-(\alpha{N}+ia)^n+(-\alpha{N}-ia)^n]}}{n!}{\int^{\infty}_{0}{\frac{dx\cdot{x^n}}{(e^{x}-1)}}}} \label{eqn1858}
\end{equation}
We recognize the Riemann zeta function and proceed to 
\begin{equation}
\sum_{k=1}^{N}{\frac{a}{(k^2+a^2)}}=\frac{1}{2i}\sum_{n=1}^{\infty}{{\alpha^{n-1}\zeta(n+1)}{[(ia)^n-(-ia)^n-(\alpha{N}+ia)^n+(-\alpha{N}-ia)^n]}} \label{eqn1859}
\end{equation} 
However, this expression will diverge with most combinations of the parameter values. The reason is that the interchange of integration and infinite summation above is not legal. This example points out when to be cautious while applying the method.
\subsubsection{Basic Trigonometric Finite Series} 
In the following we solve the $sin(x)$ finite series by using the equations ($\ref{eqn410}$) and ($\ref{eqn415}$). Here the parameters obey $0<\theta< 2\pi$
\begin{equation}
\sum_{k=1}^{N}{sin(\theta{k})}=\sum_{k=1}^{N}{g(k)} \  \label{eqn2850}
\end{equation}
The inverse Laplace transform of $g(k)$ is
\begin{equation}
G(t)=\frac{[\delta{(t+i\theta)}-\delta{(t-i\theta)}]}{2i} \label{eqn2852}
\end{equation}
We get after evaluating the integral
\begin{equation}
\sum_{k=1}^{N}{sin(\theta{k})}=\frac{cot(\frac{\theta}{2})}{2}-\frac{cos(\theta{(N+\frac{1}{2})})}{2sin(\frac{\theta}{2})} \label{eqn2855}
\end{equation} 
It will be interesting to process similarly the companion $cos(x)$ finite series
\begin{equation}
\sum_{k=1}^{N}{cos(\theta{k})}=\sum_{k=1}^{N}{g(k)}  \  \label{eqn2860}
\end{equation}
The inverse Laplace transform of $g(k)$ is
\begin{equation}
G(t)=\frac{[\delta{(t+i\theta)}+\delta{(t-i\theta)}]}{2} \label{eqn2872}
\end{equation}
and we get
\begin{equation}
\sum_{k=1}^{N}{cos(\theta{k})}=\frac{-1}{2}+\frac{sin(\theta{(N+\frac{1}{2})})}{2sin(\frac{\theta}{2})} \label{eqn2875}
\end{equation}
Both of these cases can be verified by elementary means.
\subsubsection{The Exponential Cosine Finite Series} 
In the following we tackle the negative exponentially modulated $cos(x)$ finite series in the same way as above, with $0<\theta< 2\pi$
\begin{equation}
\sum_{k=1}^{N}{e^{-\beta{k}}cos(\theta{k})}=\sum_{k=1}^{N}{g(k)}   \label{eqn8850}
\end{equation}
The inverse Laplace transform of $g(k)$ is
\begin{equation}
G(t)=\frac{[\delta{(t-(\beta-i\theta))}+\delta{(t-(\beta+i\theta))}]}{2} \label{eqn8852}
\end{equation}
We get after evaluating the integral
\begin{equation}
\sum_{k=1}^{N}{e^{-\beta{k}}cos(\theta{k})}=\frac{1}{2}[\frac{(1-e^{-N(\beta-i\theta)})}{(e^{(\beta-i\theta)}-1)}+\frac{(1-e^{-N(\beta+i\theta)})}{(e^{(\beta+i\theta)}-1)}] \label{eqn8855}
\end{equation}
This is the correct result.
\subsubsection{The Diverging Cosine Finite Series} 
Next we have the index-modulated $cos(x)$ finite series, with $0<\theta< 2\pi$, being a more complicated case.
\begin{equation}
\sum_{k=1}^{N}{{k}cos(\theta{k})}=\sum_{k=1}^{N}{g(k)}   \label{eqn18050}
\end{equation}
Our $g(k)$ is in the Laplace sense a function $f(k)$ multiplied with the argument delivering a derivative from the inverse transform. We take $f(k)=cos(\theta{k})$ and we get the inverse Laplace transform of $g(k)$ is
\begin{equation}
G(t)=\frac{d}{dt}\frac{[\delta{(t+i\theta)}+\delta{(t-i\theta)}]}{2} \label{eqn18152}
\end{equation}
We get after evaluating the integral by partial integration
\begin{equation}
\sum_{k=1}^{N}{k}cos(\theta{k})=\frac{-1}{2}[\frac{Ne^{Ni\theta}}{e^{\frac{-i\theta}{2}}(e^{\frac{-i\theta}{2}}-e^{\frac{i\theta}{2}})}+\frac{Ne^{-Ni\theta}}{e^{\frac{i\theta}{2}}(e^{\frac{i\theta}{2}}-e^{\frac{-i\theta}{2}})}] \label{eqn18255}
\end{equation}
We can process the rest in an elementary way and get the correct result
\begin{equation}
\sum_{k=1}^{N}{{k}cos(\theta{k})}=\frac{1}{2}[\frac{{N}sin(\frac
{\theta}{2})sin(\theta(N+\frac{1}{2}))-sin^2(\frac{N\theta}{2})}{sin^2(\frac{\theta}{2})}]   \label{eqn18550}
\end{equation}

\subsection{The Extended Form}
We can parametrize equation (\ref{eqn364}) with $\alpha{\in{C}}$ and $\Re(\alpha)>0$ as follows
\begin{equation}
\sum_{k=1}^{N}{g(\alpha{k})}=\int^{\infty}_{0}{\frac{dt{G(t)({1-e^{-\alpha{tN}}})}}{e^{\alpha{t}}-1}}=f(\alpha) \label{eqn380}
\end{equation}
The derivation is similar to equation (\ref{eqn364}). Since $f(\alpha)$ is a function of the parameter, we could expect it to be the result of a Laplace transform. However, its existence is not self-evident. An inverse Laplace transform would give
\begin{equation}
F(x)=\textsl{L}^{-1}_{\alpha}[f(\alpha)],x \label{eqn390}
\end{equation}
The left-hand side of equation (\ref{eqn380}) is subjected to the same transform getting
\begin{equation}
F(x)=\sum_{k=1}^{N}{\frac{G(\frac{x}{k})}{k}} \label{eqn405}
\end{equation}
These equations are equivalent, being either Laplace transforms or inverse transforms of each other.  The structure is as follows
\begin{equation}
\int^{\infty}_{0}{\frac{dt{G(t)({1-e^{-\alpha{tN}}})}}{e^{\alpha{t}}-1}}\  \  =\  \  \sum_{k=1}^{N}{g(\alpha{k})} \  \  =\  \  f(\alpha) \label{eqn410}
\end{equation}
\begin{equation}
       {\textsl{L}_{x},{\alpha}} \Uparrow \  \  \   \Downarrow{\textsl{L}^{-1}_{\alpha},{x}}   \nonumber
\end{equation}
\begin{equation}
\sum_{k=1}^{N}{\frac{G(\frac{x}{k})}{k}}\  \  =\  \  F(x), \ \ \ N<{\infty}\ \ \label{eqn415}
\end{equation}
Any of these expressions can be converted to any of the others. The functions $F(x)$ and $f(\alpha)$ are solutions to the finite series above. The equations do not imply that any function could be expanded as a finite series by using these equations. Appendix A shows variants of equations (\ref{eqn410}) and (\ref{eqn415}) for some special cases. 

The equations for finite series are interesting since convergence is not required and the number of terms can be brought down to one. This allows almost any arbitrary function of the index to be applied. The result may become a closed-form function or another series. 
\subsection{The Extended Method}
We may have a finite series of type A in a parametrized form 
\begin{equation}
\sum_{k=1}^{N}{g(\alpha{k})} \label{eqn428}
\end{equation}
or we might have a finite series similar to type B
\begin{equation}
\sum_{k=1}^{N}{\frac{G(\frac{x}{k})}{k}} \label{eqn430}
\end{equation}
It is possible to proceed in different ways to solve a finite series:

	\begin{itemize}
  \item Type A,  equation (\ref{eqn428}), integral
	\begin{itemize}
	\item 1. starting from equation (\ref{eqn410})
	\item  2. determine the inverse Laplace transform $G(t)$
	\item  3. solve the integral in equation (\ref{eqn410}) and get $f(\alpha)$
	\item  4. the parameter $\alpha$ is set to $1$. If only equation (\ref{eqn410}) is applied, the parameter $\alpha$ is not necessary. The solution is $f(1)$
  \end{itemize}

  \item Type A,  equation (\ref{eqn428}), via type B finite series
	\begin{itemize}
  \item 1. starting from equation (\ref{eqn410}) 
  \item 2. calculate the inverse Laplace transform $G(t)$
  \item 3. solve the new finite series (\ref{eqn415}) getting $F(x)$ 
  \item 4. $F(x)$ is Laplace transformed to get $f(\alpha)$. The solution is $f(1)$	
  \end{itemize}

  \item Type B finite series equation (\ref{eqn430}), via type A finite series
	\begin{itemize}
  \item 1. starting from equation (\ref{eqn415}) 
  \item 2. generate the Laplace transform $g(k)$ from $G(t)$
	\item 3. parametrize it to $g(\alpha{k})$ 
	\item 4. solve the finite series in equation (\ref{eqn410})
	\item 5. inverse transform the resulting $f(\alpha)$ to get $F(x)$, the solution. 
	\item 6. Set $x$ to some final value depending on the fitting in equation (\ref{eqn415}).
  \end{itemize}

  \item Type B finite series equation (\ref{eqn430}), via the integral and inverse transform
	\begin{itemize}
  \item 1. starting from equation (\ref{eqn415}) 
  \item 2. use the $G(t)$ to solve the integral in equation (\ref{eqn410}) to get $f(\alpha)$
	\item 3. inverse transform $f(\alpha)$ to get $F(x)$. 
	\item 4. Set $x$ to some final value depending on the fitting in equation (\ref{eqn415}).
  \end{itemize}
  \end{itemize}
Verification of the result should be the final task. Success of this method is dependent on existence and eventual finding of the Laplace transforms and inverse transforms necessary, on solving the new finite series and/or the resulting integral. 

In some cases there may appear infinite series. Interchanging integration and infinite summations must be studied carefully since they may lead to invalid results or divergence. As is well known, it is required that the infinite series must be uniformly convergent for interchanging integration and summation of its terms.
\section{The Fourier Transform Method}
By using the Fourier transform pair below we can generate an analogous method as in the preceding section. 
\begin{equation}
F_{x}[g(x)]_{,\alpha}=G(\alpha)=\int^{\infty}_{-\infty}{dx{g(x)e^{-i\alpha{x}}}} \label{eqn480}
\end{equation}
\begin{equation}
F^{-1}_{\alpha}[G(\alpha)]_{,x}=g(x)=\frac{1}{2\pi}\int^{\infty}_{-\infty}{d\alpha{G(\alpha)e^{i\alpha{x}}}} \label{eqn490}
\end{equation}
It can be applied with the following traditional conditions (convergence and Dirichlet)
\begin{itemize}
	\item 1. the integral $\int^{\infty}_{-\infty}{dx{\left|g(x)\right|}}$ converges
	\item 2. $g(x)$ and $\frac{dg(x)}{dx}$ are piecewise continuous in every finite interval $-L < x < L$
	\item 3. at discontinuities $g(x)$ is replaced with $\frac{1}{2}[g(x+0)+g(x-0)]$
	\item 4. $g(x)$ is defined and single-valued except possibly at a finite number of points in $(-L, L)$
	\item 5. $g(x)$ is periodic outside $(-L, L)$ with a period of $2L$
\end{itemize}
We follow the track as above by substituting the transform and use
\begin{equation}
\sum_{k=1}^{N}e^{i\alpha{k}}=\frac{1-e^{i{N\alpha}}}{e^{-i\alpha}-1} \label{eqn500}
\end{equation}
\begin{equation}
\sum_{k=1}^{N}{g(k)}=\frac{1}{2\pi}\sum_{k=1}^{N}{\int^{\infty}_{-\infty}{d\alpha{G(\alpha)e^{i\alpha{k}}}}} \label{eqn501}
\end{equation}
\begin{equation}
=\frac{1}{2\pi}\int^{\infty}_{-\infty}{\frac{d\alpha{G(\alpha)(1-e^{i\alpha{N}})}}{e^{-i\alpha}-1}} \label{eqn510}
\end{equation}
Interchanging the integration and summation is legal here. We can replace the exponential functions with $sin(x)$ functions and get
\begin{equation}
\sum_{k=1}^{N}{g(k)}=\frac{1}{2\pi}\int^{\infty}_{-\infty}{\frac{d\alpha{G(\alpha)sin(\frac{\alpha{N}}{2}})}{sin(\frac{\alpha}{2})}} \label{eqn550}
\end{equation}
An analogous generalization, as done earlier in equations (\ref{eqn410}) and (\ref{eqn415}), can be developed for this method. Presently, use of this method is limited due to rather short tables for Fourier transforms. An extensive work is required to generate proper tables to aid using this method. 

\section{A Simplified Approach}
One may, while applying the equation (\ref{eqn364}), meet obstacles. They may be due to repulsive-looking integrals and other complexities. We can handle the problem differently. We continue processing it as follows
\begin{equation}
\sum_{k=1}^{N}{g(k)}=\int^{\infty}_{0}{\frac{dt{G(t)({e^{-t}-e^{-t{(N+1)}}})}}{1-e^{-t}}} \label{eqn1424}
\end{equation}
Next, we use the binomial series expansion for the denominator to obtain
\begin{equation}
=\sum_{n=0}^{\infty}{\int^{\infty}_{0}{dt{G(t)({e^{-t(n+1)}-e^{-t{(N+n+1)}}})}}} \label{eqn1426}
\end{equation}
We recognize this expression as a pair of Laplace transforms of the $G(t)$. $g(k)$ must have values for the index beyond $N$ to infinity, a mapping for each index value. Therefore we get
\begin{equation}
\sum_{k=1}^{N}{g(k)}=\sum_{n=0}^{\infty}{[g(n+1)-g(N+n+1)]}\label{eqn1428}
\end{equation}
\begin{equation}
\sum_{k=1}^{N}{g(k)}=\sum_{k=1}^{\infty}{[g(k)-g(N+k)]}\label{eqn1430}
\end{equation}
The result is universal. This is actually a rather trivial result but offers another way of solving or transforming a finite series in those cases where a corresponding infinite series does exist. The infinite series for $g(k)$ does not need to converge since we are subtracting term-by term and the separate sums may diverge. We still govern the $N$ freely.
\subsection{Example - An Inverse Power Series} 
The finite series of the type of Riemann zeta function is handled with equation (\ref{eqn1430}). 
\begin{equation}
\sum_{k=1}^{N}{\frac{1}{k^{s}}}=\sum_{k=1}^{N}{g(k,s)} \  \label{eqn3000}
\end{equation}
Since for $\Re{(s)}>1$
\begin{equation}
\sum_{k=1}^{\infty}{\frac{1}{k^{s}}}=\zeta(s) \  \label{eqn3010}
\end{equation}
and the corresponding Hurwitz zeta function is
\begin{equation}
\sum_{k=1}^{\infty}{\frac{1}{(k+N)^{s}}}=\zeta(s,N)-\frac{1}{N^s} \  \label{eqn3020}
\end{equation}
we obtain a simple result
\begin{equation}
\sum_{k=1}^{N}{\frac{1}{k^{s}}}=\zeta(s)-\zeta(s,N)+\frac{1}{N^s} \  \label{eqn3030}
\end{equation}
\section{Discussion}
We have presented methods for generating transformations and solutions to general finite series. Various types of functions $g(k)$ of the index can be handled. The number $N$ of the terms in the finite series can be running from one up to some high number. However, if infinite series are required or the limit $N\rightarrow{\infty}$ needs to be taken, the method may fail or the equations become unsurmountable. In such cases the infinite series method in \cite{Stenlund20142} may produce correct results.

The main results of this study are equations (\ref{eqn364}), (\ref{eqn410}), (\ref{eqn415}), (\ref{eqn550}) and (\ref{eqn1430}). We have a simpler form of solution for finite series in equation (\ref{eqn364}) and a generalized form in equations (\ref{eqn410}) and (\ref{eqn415}). A parameter $\alpha$ can be added to the functions as a multiplier of the index. In the simple form of the equations, it is not required. The results bear similarity to the infinite series method but the number of terms $N$ in finite series will bring in extra terms to the equations. The requirements on the functions $g(k)$ and $G(t)$ are not as severe as with infinite series. We display a variant of the method realized with the Fourier transform, equation (\ref{eqn550}). Its use is more limited due to the small number of tabularized transforms. We present a simplified method in equation (\ref{eqn1430}) which is very powerful in some cases.

It is important for a new solution to be verified. Success of the method is dependent on finding the necessary Laplace transforms or inverse transforms. Solving the integral in equation (\ref{eqn410}) may be complicated too and easily leads to an infinite series of Riemann zetas if series expansions are used as an aid. The method may fail if caution is not taken at every step, especially while interchanging integration and infinite summation.

These methods are important for finding closed-form expressions representing the series. New interesting identities can be generated since the number of terms is limited.

\newpage
\appendix
\section{Appendix. Various Particular Forms of Equivalence}
We can generate different equations modifying equation (\ref{eqn410}) and (\ref{eqn415}) to help solving or transforming finite series. Proofs are similar to those before. We have added subscripts to the functions $f(\alpha)$ and $F(x)$ to remind that they differ between cases. The formulas are given names referring to their expected crude behavior as a finite series. Obviously, it is straightforward to generate more of these equations but they seem not to be so useful due to their complexity.

\subsection{Alternating Finite Series}
These equations are valid \textbf{only for an even N}.
\begin{equation}
\int^{\infty}_{0}{\frac{dt{G(t)(1-e^{-\alpha{tN}})}}{e^{\alpha{t}}+1}}\  \  =\  \  \sum_{k=1}^{N}{(-1)^{k+1}g(\alpha{k})} \  \  =\  \  f_{1}(\alpha) \label{eqn5000}
\end{equation}
\begin{equation}
       {\textsl{L}_{x},{\alpha}} \Uparrow \  \  \   \Downarrow{\textsl{L}^{-1}_{\alpha},{x}}  \nonumber
\end{equation}
\begin{equation}
\sum_{k=1}^{N}{\frac{(-1)^{k+1}G(\frac{x}{k})}{k}}\  \  =\  \  F_{1}(x) \label{eqn5010}
\end{equation}

\subsection{Shifted Finite Series}
\begin{equation}
\int^{\infty}_{0}{\frac{dt{G(t)(1-e^{-\alpha{tN}})}e^{-\beta{t}}}{e^{\alpha{t}}-1}}\  \  =\  \  \sum_{k=1}^{N}{g(\alpha{k}+\beta)} \  \  =\  \  f_{2}(\alpha) \label{eqn5020}
\end{equation}
\begin{equation}
       {\textsl{L}_{x},{\alpha}} \Uparrow \  \  \   \Downarrow{\textsl{L}^{-1}_{\alpha},{x}}  \nonumber
\end{equation}
\begin{equation}
\sum_{k=1}^{N}{\frac{e^{\frac{\beta{x}}{k}}G(\frac{x}{k})}{k}}\  \  =\  \  F_{2}(x) \label{eqn5030}
\end{equation}

\subsection{Shifted Alternating Finite Series}
These equations are valid \textbf{only for an even N}.
\begin{equation}
\int^{\infty}_{0}{\frac{dt{G(t)(1-e^{-\alpha{tN}})}e^{-\beta{t}}}{e^{\alpha{t}}+1}}\  \  =\  \  \sum_{k=1}^{N}{(-1)^{k+1}g(\alpha{k}+\beta)} \  \  =\  \  f_{3}(\alpha) \label{eqn5040}
\end{equation}
\begin{equation}
       {\textsl{L}_{x},{\alpha}} \Uparrow \  \  \   \Downarrow{\textsl{L}^{-1}_{\alpha},{x}}  \nonumber
\end{equation}
\begin{equation}
\sum_{k=1}^{N}{\frac{e^{\frac{\beta{x}}{k}}(-1)^{k+1}G(\frac{x}{k})}{k}}\  \  =\  \  F_{3}(x) \label{eqn5050}
\end{equation}

\subsection{Exponential Factor Finite Series}
With $\Re(\beta)>0$
\begin{equation}
\int^{\infty}_{0}{\frac{dt{G(t)}(1-e^{-N(\beta+\alpha{t})})}{e^{\alpha{t}+\beta}-1}}\  \  =\  \  \sum_{k=1}^{N}{e^{-\beta{k}}g(\alpha{k})} \  \  =\  \  f_{4}(\alpha) \label{eqn6000}
\end{equation}
\begin{equation}
       {\textsl{L}_{x},{\alpha}} \Uparrow \  \  \   \Downarrow{\textsl{L}^{-1}_{\alpha},{x}}   \nonumber
\end{equation}
\begin{equation}
\sum_{k=1}^{N}{\frac{e^{-\beta{k}}G(\frac{x}{k})}{k}}\  \  =\  \  F_{4}(x) \label{eqn6010}
\end{equation}

\subsection{Exponential Factor Alternating Finite Series}
These equations are valid \textbf{only for an even N}. With $\Re(\beta)>0$
\begin{equation}
\int^{\infty}_{0}{\frac{dt{G(t)}(1-e^{-N(\beta+\alpha{t})})}{e^{\alpha{t}+\beta}+1}}\  \  =\  \  \sum_{k=1}^{N}{(-1)^{k+1}{e^{-\beta{k}}g(\alpha{k})}} \  \  =\  \  f_{5}(\alpha) \label{eqn6020}
\end{equation}
\begin{equation}
       {\textsl{L}_{x},{\alpha}} \Uparrow \  \  \   \Downarrow{\textsl{L}^{-1}_{\alpha},{x}}   \nonumber
\end{equation}
\begin{equation}
\sum_{k=1}^{N}{\frac{(-1)^{k+1}e^{-\beta{k}}G(\frac{x}{k})}{k}}\  \  =\  \  F_{5}(x) \label{eqn6030}
\end{equation}

\end{document}